\documentclass[11pt,a4paper]{article}
\usepackage{t1enc}
\usepackage[latin1]{inputenc}
\usepackage[german]{babel}
\usepackage{harvard}
\usepackage{graphicx}

\newtheorem{Sachverhalt}{Sachverhalt}

\newcommand{\beq}{\begin{equation}}
\newcommand{\eeq}{\end{equation}}

\begin{document}
\thispagestyle{empty}
\begin{center}
\subsection*{C.F. Gau\ss{}' Pr\"azisionsmessungen terrestrischer Dreiecke und
seine Überlegungen zur empirischen Fundierung der Geometrie in den 1820er
Jahren  }
{\em Erhard Scholz, Wuppertal}\footnote{scholz@math.uni-wuppertal.de}
\end{center}
\begin{abstract}
In the historical literature there has been an extended discussion on the question, whether  the report of Sartorius von Waltershausen about C. F. Gauss checking the largest  triangle of the geodetical measurement campaign in the kingdom of Hannover as  a kind of ``test'' for the Euclididean nature of  physical space can be taken seriously or not.  Among others, it was argued that it even  was  logically  impossible for Gauss  to do so in the early 1820s (i.e. in particular, before J. Bolyai's, N.I. Lobachevsy's, or even B. Riemann's works). This article shows, in which sense Gauss's methodology of  curvature of surfaces, although logically developed in all clarity only for surfaces embedded in Euclidean 3-space, could very well be used  already in the early 1820   to investigate the above mentioned question  in the sense  of physical geometry. Although we do not have a definitive  proof of the respective calculations by  Gauss's own hand, the latter's  account of the situation to contemporaries (correspondents, friend and students) changed clearly between the early and the late 1820s. That speaks very much in favour of a basic correctness of Sartorius von Waltershausen's report on this topic.

\end{abstract}

\subsubsection*{1. Zum Stand der historischen Diskussion}
In dem von unserem Jubilar  herausgegebenen Sammelband zum 200. Geburtstag
von Carl Friedrich Gau\ss{} \cite{Schneider:Gauss} stellte H. Gericke die Gau\ss
schen Überlegungen zur Rolle der  ``charakteristischen Konstante'' $C >
0$ der nichteuklidischen Geometrie dar. In differentialgeometrischer Sprache
wird durch sie die Krümmung  charakterisiert, $\kappa = - C$, und die
Winkelsumme eines Dreiecks weicht bei gegebener Konstante $C$  um einen vom
Fl\"acheninhalt $F$ abh\"angigen Betrag $\epsilon $ vom euklidischen Wert ab,
\[\epsilon   = 180^{\circ} - ( \alpha +\beta +\gamma )   =   C \cdot F . \]

Herr Gericke zitierte in diesem Kontext aus einem Brief vom 8. 11. 1824 des
damaligen Jubilars, Gau\ss{}, an Taurinus : ``W\"are die Nicht-Euklidische
Geometrie die wahre, und jene Constante in einigem Verh\"altnisse zu solchen
Grö\ss{}en, die im Bereich unserer Messungen auf der Erde oder am Himmel
liegen, so lie\ss{}e sie sich  a posteriori ausmitteln''
\cite[VIII, 186-188]{Gauss:Werke}.
H. Gericke  gab zu  dieser treffenden Bemerkung über die Möglichkeit einer
empirischen  ``Ausmittelung'' der physikalischen Raumkrümmung den Kommentar:
\begin{quote}
Man könnte z.B. die Winkelsumme in einem Dreieck messen. Da der Defekt
[oben als $\epsilon  $ bezeichnet, E.S.] dem Fl\"acheninhalt proportional ist,
mu\ss{} man ein möglichst gro\ss{}es Dreieck w\"ahlen. Gau\ss{} w\"ahlte bekanntlich
das Dreieck Brocken -- Inselsberg -- Hoher Hagen, dessen grö\ss te Seite
ca. 105 km lang ist, und das Ergebnis war  bekanntlich: Der Unterschied der
Winkelsumme von $180^{\circ}$ war jedenfalls kleiner als der Me\ss
fehler.\footnote{Gau\ss{} hatte im Jahre 1823, also im Vorjahr des zitierten
Briefes an Taurinus, das erw\"ahnte Dreieck mit einer bis dahin nicht
erreichten Pr\"azision ausgemessen und eine Abweichung der Winkelsumme von
etwa 0,6'' erhalten; weitere Details siehe unten. Anm. E.S.}
  Also: entweder ist die euklidische Geometrie wahr oder das Dreieck war
immer noch zu klein. Die um eine Entscheidung befragte Natur hatte die
Antwort verweigert. \cite[129]{Gericke:Gauss}
\end{quote}

Dieser Kommentar gab   die bis Anfang der 1970er Jahre  weitgehend
unumstrittene Auffassung über die (geometrisch--) grundlagentheoretische
Auswertung der geod\"atischen Messungen durch Gau\ss{} wieder.
Um diese Zeit wurde jedoch unter einigen Wissenschaftshistorikern die
Verl\"asslichkeit der Überlieferung in diesem Punkt  angezweifelt. Gau\ss{}
hatte sich in den uns überlieferten schriftlichen Äu\ss{}erungen nie
definitiv zu dieser speziellen Frage  einer raumstrukturellen
Auswertung seiner geod\"atischen Messungen ge\"au\ss{}ert. Die angegebene
Auffassung basierte also neben einer naheliegenden historischen
Extrapolation aus anderen Quellen, wie der oben zitierten Korrespondenz, auf
den Anmerkungen des Zeitzeugen und Freundes von C.F. Gau\ss{}, Sartorius von
Waltershausen, der in seinem Band {\em Gau\ss{} zum Ged\"achtnis} kurz nach
dessem Tod berichtet hatte:
\begin{quote}
Die Geometrie betrachtete Gauss nur als ein consequentes Geb\"aude nachdem die
Parallelentheorie als Axiom an der Spitze zugegeben sei: er sei indess zur
Überzeugung gelangt, dass dieser Satz nicht bewiesen werden könne, doch
wisse man aus der Erfahrung z.B. aus den Winkeln des Dreiecks Brocken,
Hohenhagen, Inselsberg, dass er n\"aherungsweise richtig sei. 
Wolle man dagegen das genannte Axiom nicht zugeben, so folge daraus eine andere ganz selbstst\"andige Geometrie, die er gelegentlich ein Mal verfolgt und mit dem Namen Antieuklidische Geometrie bezeichnet habe. \cite[81]{Sartorius}
\end{quote}
An anderer Stelle p\"azisierte Sartorius die vage, für sich genommen nicht
sonderlich aufschlussreiche Aussage ``n\"aherungsweise richtig'' durch eine
Anmerkung zu Gau\ss{}' neu entwickeltem geod\"atischen
Pr\"azisionsmessinstrument, dem Heliotrop:
\begin{quote}
Das Heliotrop fand sogleich bei der Hannöverschen Triangulation seine volle
Anwendung und das grosse Dreieck, vielleicht noch das grösste, welches
gemessen worden ist, n\"amlich zwischen dem Brocken, dem Inselsberg und dem
Hohenhagen, wurde mit Hülfe desselben so  genau gemessen, dass die Summe der
drei Winkel nur etwa um zwei Zehntheile einer Secunde sich von zwei Rechten
entfernt. \cite[53]{Sartorius}
\end{quote}
Diese Pr\"azisierung,   ungef\"ahr ``zwei Zehntheile einer Secunde'', ist ein
deutlicher Hinweis auf ein von Gau\ss{} in mündlichen Gespr\"achen
mitgeteiltes Pr\"azisionsma\ss{} bei der Auswertung seiner Messungen. Es w\"are
sehr befremdlich zu unterstellen, dass der Berichterstatter  eine solche
Angabe einer Fehlerschranke erfunden und seinem Bericht zugefügt h\"atte. Es
bleibt aber zu kl\"aren, um welche Art der Fehlerschranke   es
Gau\ss{} bei der von Sartorius zitieren Bemerkung ging oder gehen konnte.
Wir werden darauf zurückzukommen haben.

Im Jahre 1972 stellte der damals junge Physikhistoriker A. Miller  in einem
Artikel in {\em Isis} die bis dahin weitgehend unangefochtene Darstellung
des Sartorius von Waltershausen in Frage \cite{Miller:Myth}. Dieser Artikel
erzielte eine durchschlagende Wirkung.  In einem anderen, von K. Reich
verfassten, Beitrag zum selben  Sammelband, in dem H. Gericke auf die
``wohlbekannte'' Geschichte der Auswertung verwies, fügte die Autorin der
Wiedergabe der ersten der beiden oben wiedergegebenen Zitatstellen bei
Sartorius folgende kritische Bemerkung an:
\begin{quote}
In der Sekund\"arliteratur sind oftmals Bemerkungen der Art zu finden, da\ss{}
Gau\ss{} bei der Vermessung dieses Dreiecks die wahre Struktur des Raumes,
sei er euklidisch oder nichteuklidisch in Erfahrung bringen wollte. Da
Gau\ss{} selbst diese Verbindung nicht ausdrücklich hergestellt hat, ist
diese Interpretation vor einigen Jahren sowohl angezweifelt als auch
best\"atigt und wieder angezweifelt worden.
\cite[104f.]{Reich:Gauss_Schneider}
\end{quote}
Damit wiedersprach sie offen, wenn auch zurückhaltend  der Darstellung von
Herrn Gericke im selben Band. Zwar enthielt sie sich einer  eigenen
Stellungnahme zu diesem Problem, verwies in einer Anmerkung aber auf die
damals vorliegende Literatur.\footnote{K. Reich gab folgende
Literaturstellen an: \cite{Miller:Myth} (``angezweifelt''),
\cite{Goe:Myth,Waerden:Myth} (``best\"atigt'') und  \cite{Miller:Myth2} (``wieder angezweifelt''). }
Bei einer anderen Gelegenheit habe ich selber schon  kurz erkl\"art, warum ich
die gegen die Darstellung des Sartorius von Waltershausen   vorgetragenen
Zweifel  nicht teile, also in K. Reichs Formulierung die erw\"ahnte Geschichte
als ``wieder best\"atigt'' ansehe  \cite[643f.]{Scholz:Amphora}.  Die
Diskussion unter Wissenschaftshistorikern wurde also fortgeführt,
jedoch, wie es scheint,  mit  weiterhin offenem  Ergebnis.  Die 
detaillierteste, höchst  sachkundige Darstellung findet sich bis heute bei
\cite{Breitenberger:Gauss}; sie beschr\"ankt sich allerdings darauf zu belegen,   ``was Gau\ss{} tats\"achlich getan hat und was man ihm ungestraft nachsagen darf''.\footnote{Persönliche Mitteilung E. Breitenberger an den Autor, 18. 02. 2003.}
Dabei entsteht  ein --- aus meiner Sicht übertrieben --- skeptisches historisches
 Gesamtbild der Gau\ss{}schen Perspektive.\footnote{Ich verdanke Jeremy Gray den Hinweis
auf diese ``offene Situation'' und damit in gewissem Sinne die Anregung zu
diesen Ausführungen.  In einer Serie von Diskussionen hat er mir  klar
gemacht,   warum man bessere Gründe des Zweifels haben kann, als in dem
ursprünglichen Artikel von A. Miller angeführt.
Bei einen dieser Gespr\"ache (Tel Aviv  im Mai 2001) war unser
Jubilar anwesend und zeigte sich daran so interessiert, dass sich dieses
Thema  für diesen Sammelband wohl eignen mag.  }

Die derzeitige Situation l\"asst sich wohl am besten so zusammenfassen: Unter
Wissenschaftshistorikern wird tendenziell eher bezweifelt als akzeptiert,
dass Gau\ss{} seine Messungen wie angegeben ausgewertet h\"atte (Miller)
oder dass er es überhaupt  sinnvoll h\"atte tun können (Breitenberger u.a.).
Unter  Mathematikern und mathematischen Naturwissenschaftlern (insbesondere
Geod\"aten und Physikern) wird dies eher als gesichert angesehen als bezweifelt;
diese Sicherheit  beruht jedoch zum Teil  auf Unkenntnis der vorgebrachten
Zweifelsgründe.

Es mag also  gerechtfertigt sein, hier auf diese Frage zurückzukommen. Ich
werde dazu  die ursprünglichen Gegenargumente von A. Miller ``kritisch
würdigen'', in diesem Fall hei\ss{}t dies aufzeigen, warum sie am
historischen Sachverhalt insgesamt  ziemlich  weit vorbeigehen (Teil 2).
Dann werde ich die von E. Breitenberger vorgetragenen methodologischen
Probleme vorstellen, die beim Versuch  einer empirischen Überprüfung der
Euklidizit\"at des physikalischen Raumes durch geod\"atische Messungen auftreten
und in Rechnung zu stellen sind  (Teil 3). Schlie\ss{}lich möchte ich
 erkl\"aren, warum meiner Ansicht nach s\"amtliche der problematischen Punkte
{\em   im Rahmen der Gau\ss{} zur Verfügung stehenden Methoden als
ausreichend  absch\"atzbar} angesehen werden können, obwohl letzterem {\em
nicht alle theoretischen Voraussetzungen  in voll ausgearbeiteter Form} zur
Verfügung standen (Teil 4). Es wird sich zeigen, dass Gau\ss{} im Rahmen
seiner Fehlerabsch\"atzungen durchaus  in der Lage war, strukturelle Fragen
empirisch aussagekr\"aftig zu entscheiden (Teil 5). Damit kann ich keinen
nachvollziehbaren Grund mehr sehen, die Verl\"asslichkeit des Zeitzeugen
Sartorius von Waltershausen in Zweifel zu ziehen (Teil 6).

\subsubsection*{2. Die Etablierung eines Mythos}
A. Millers ursprüngliche  Kritik an der überlieferten Geschichte unter dem
Titel {\em The myth of  Gauss' experiment on the Euclidean nature of
physical space} enth\"alt aus rückblickender Sicht erstaunlich wenig
Argumente, wenn man die  erhebliche Wirkung des Artikels in Rechnung
stellt. Miller stellte folgende einfache, aber wie wir sehen werden,
schlecht belegte zentrale These auf:
\begin{quote}
The famous experiment is in fact a mere legend that stems from a
misunderstanding of his important paper of 1827  ``Disquisitiones generales
circa superficies curvas.'' This work constituted the first systematic study
of quadratic differential forms (of two variables) and was generalized by
Riemann, in 1854, to $n$ variables. \cite[346]{Miller:Myth}
\end{quote}

Miller bezog sich damit auf den \S 28 der  Arbeit \cite{Gauss:1828}, in dem
Gau\ss{} mittels des von ihm pr\"azisierten und verallgemeinerten
Legendre-Theorem der sph\"arischen Geometrie nachprüfte,\footnote{Das
Legendre-Theorem der sph\"arischen Geometrie besagt: Ein sph\"arisches Dreieck
mit kleinen Seiten und daher kleinem sph\"arischen Exzess hat approximativ
gleichen Fl\"acheninhalt wie das ebene Dreieck mit gleichen Seitenl\"angen.
Jeder Winkel des ebenen Dreiecks ist um ein Drittel des sph\"arischen Exzesses
kleiner als der entsprechende Winkel des sph\"arischen Dreiecks. --- Gau\ss{} gab
eine Pr\"azisierung der Approximation des Fl\"acheninhalts und verallgemeinerte
den Satz samt Winkelkorrekturen für den Fall eines ``kleinen'' Dreiecks auf
einer beliebig gekrümmten Fl\"ache. Vgl. \cite{Dombrowski:Gauss} oder \cite[639f.]{Scholz:Amphora}}
 ob die Nicht\-sph\"a\-ri\-zi\-t\"at der Erdgestalt schon bei
Grö\ss{}enordnungen der von ihm ausgemessenenen Dreiecke   quantitative
Auswirkungen in  N\"ahe der Messgenauigkeit hat und daher  ggfs. schon bei der
{\em Auswertung der Daten eines einzelnen geod\"atischen Dreiecks} in Rechnung
zu stellen ist.

Gau\ss{} verglich dabei die Winkelkorrekturen, die an einem Dreieck
$\triangle $ auf der sph\"arisch oder ellipsoidisch angenommenen Erdoberfl\"ache
angebracht werden müssen, um zu einem seitengleichen ebenen (euklidischen)
Dreieck $\triangle ^{\ast }$ überzugehen. Im sph\"arischen Fall ist der
Winkelüberschuss des Dreiecks $\triangle $ über $180^{\circ}$ gleich\-m\"a\ss{}ig
auf die drei Winkel zu verteilen, im nichtsph\"arischen Fall verschieden, in
Abh\"angigkeit von der Krümmung der Fl\"ache (Erdgestalt) in   den einzelnen
Ecken. Gau\ss{} ermitttelte, dass im Testfall des Dreiecks Brocken -
Hohehagen - Inselsberg (das ich im folgenden  kurz als das {\em  gro\ss{}e
Dreieck} oder $\triangle BHI$ bezeichnen möchte) der von ihm ermittelte
Winkelüberschuss $14'',85348$ in den Teilbetr\"agen $4'',95104, \, 4'',95113$
und $ 4'',95131$ auf die Ecken $B, H, I$ aufzuteilen ist, also erst in der
vierten Nachkommastelle einer Winkelsekunde und damit  weit unterhalb der
Messgenauigkeit (erste Nachkommastelle bei den besten Messungen)\footnote{Der mittlere
Fehler (im Sinne der empririschen Standardabweichung) war
 $3,5'' /\sqrt{n}$ bei $n$ Messungen derselben
Richtung \cite[273, 279f., 281]{Breitenberger:Gauss}} von der
Gleichverteilung abweicht.

Miller gab eine zutreffende  --- nebenbei angemerkt, vorher   keineswegs
bestrittene --- Darstellung dieses Sachverhaltes:
\begin{quote}
That Gauss concluded nothing about the non-Euclidean nature of physical
space
[in \cite[\S 28]{Gauss:1828}, E.S.] is not surprising, because the
mathematical theory which he developed in this paper was not at all
concerned with non-Euclidean geometry. What Gauss was seeking was a
generalization of Legendre's result to a doubly curved surface in order to
determine whether the earth 's double curvature had any effect on his
geodetic data.  \cite[348]{Miller:Myth}
\end{quote}
Er fuhr fort:
\begin{quote}
Thus, Gauss did not have to use his geodetic data to determine whether the
space of our experience is curved. He was studying a curved surface of known
curvature --- the earth. Additional proof that a misunderstanding of Gauss'
1827 paper led to the legend of his ``experiment'' is that of the many
triangles that he surveyed and reported on, the previously mentioned
triangle is the only one taken into account in the aforementioned paper. It
therefore appears that what up to now has often been presented as Gauss'
experiments on the nature of the space in which we live is simpy a myth.
(ibid.)
\end{quote}
Miller widerlegte hier überraschenderweise ein Missverst\"andnis,
(``misunderstanding of Gauss' 1827 paper''), das in der informierten
mathematikhistorischen Literatur vor ihm  keine  Rolle
spielte.
 Weder Sartorius von Waltershausen, noch einer der daran anknüpfenden
informierten Autoren hatte vor A. Miller den \S 28 der {\em Disquisitiones
generales} mit der  Nicht\-eukli\-dizi\-t\"ats\-prü\-fung durch  Gau\ss{} in Beziehung
gebracht. Das nach Sartorius  von Gau\ss{} erw\"ahnte Genauigkeitsma\ss{}
 (``um zwei Zehntheile  einer Secunde von zwei Rechten entfernt'')  widersprach bei genauem Hinsehen ziemlich in allem einer solchen  Interpretation des \S 28: relevante
quantitative Abweichung in der Grö\ss{}enordnung
$10^{-4}$ Bogensekunden statt $10^{-1}$, eine Abweichung von zwei Rechten
stand hier  nicht zur Debatte, vielmehr ging es hier um eine ``Umrechnung''
auf zwei Rechte im oben angemerkten Sinne der Bestimmung eines
seitengleichen ebenen Dreieckes $\triangle ^{\ast }$.\footnote{Dieses Argument führte schon van der Waerden in seiner Kritik von 1974 an.} Das Hauptargument  A.
Millers richtete sich also  lediglich auf ein   {\em von ihm selbst  in
die Historiographie der Mathematik hineingetrages Mi\ss{}verst\"andnis.} Der
Hinweis auf das ``gro\ss{}e'' Dreieck in Millers   ``additional proof'' geht
völlig an der naheliegenden Tatsache vorbei, dass  Gau\ss{} dieses
möglicherweise aus guten (und leicht nachvollziehbaren)  Gründen   gerne für
Kontrollüberlegungen {\em verschiedener Art} heranzog --- schlie\ss{}lich hatte
er besondere Sorgfalt für eine möglichst genaue und so weit wie möglich vom Rest
des Netzes unabh\"angige  Messung dieses Dreiecks  aufgewendet.\footnote{Siehe
dazu \cite{Gerardy:Gauss}.}

Miller erhob dabei gar nicht den Anspruch, den Bericht des  Zeitzeugen
Sartorius von Waltershausen zu kritisieren. Anscheinend war ihm  
zur Zeit der Abfassung seines
 ``Myth''-Artikels  dieser Bericht sogar g\"anzlich unbekannt. Sonst w\"are
schwer zu verstehen und als die Leser  absichtlich in die  Irre führend anzusehen,
dass im gesamten Artikel   \cite{Miller:Myth} Sartorius nicht erw\"ahnt wird,
obwohl dieser der zentrale Bezugspunkt und ``heimliche Gegner''  des
Millerschen  Widerlegungsversuches war.\footnote{In seiner Antwort an die Kritiker verteidigte sich Miller dann mit dem Versuch einer Abwertung des Zeitzeugenberichtes mit m.E. fadenscheinigen Sekund\"arargumenten vom Typ: auch Ernst Mach habe Sartorius nicht ernst genommen usw. \cite{Miller:Myth2}.}   Anstelle eines  Bezugspunktes fand
Miller jedoch nur eine Lücke vor, die er durch eine kühne  Hypothese  über
den Ursprung der von ihm zum ``Mythos'' erhobenen Geschichte  zu
schlie\ss{}en versuchte:
\begin{quote}
In conclusion, a reasonable conjecture as to the origin of this myth is as
follows. The question of whether the physical space is curved or not took on
new meaning after Einstein's general theory of relativity (1916), which
utilized Riemannian geometry. Consequently, the myth may have arisen as a
result of extrapolating back from Einstein's work to Riemann's 1854 work to
Gauss' 1827 paper (without reading it, of course), keeping in mind that
Gauss's theorems were applied to the largest triangle in his geodetic
survey. \ldots \cite[348]{Miller:Myth}
\end{quote}
Das war  ein mehr als  eigenwilliger Versuch einer Neukonstruktion der
Überlieferungsgeschichte, der die gesamte vor-relativistische Debatte um die
``Natur'' des Raumes im sp\"aten 19. und frühen 20. Jahrhundert komplett
ignorierte.  Es w\"are deplaziert, die hier vorgestellte gegenstandslose
Hypothese  ausdrücklich ``widerlegen'' zu wollen; ihr Autor wird sie in der
folgenden Debatte schrittweise wieder aufgegeben haben, nachdem (oder soll ich sagen ``falls''?) er von der
vor-relativistischen
Diskussion über die empirischen Bezüge der nichteuklidischen Geometrie
erfahren hat.
In ganz anderem Sinnne  mag diese Passage allerdings aufschlussreich sein,
n\"amlich als ein  Indikator des  Hintergrundes, aus dem heraus dieser Artikel
geschrieben wurde und seine Wirksamkeit entfalten konnte. Zun\"achst einmal
muss es ja  r\"atselhaft erscheinen, wie es kommen konnte, dass ein so
 uninformierter und schlecht recherchierter Artikel als
Standard\-referenz auch für die sp\"atere, besser begründete  Kritik an der
Überlieferungsgeschichte der Gau\ss{}schen Euklidizit\"atsüberprüfung anhand
des ``gro\ss{}en Dreiecks'' dienen konnte.

Ich möchte hier nur kurz zusammenfassen: Die im Bericht des Sartorius von
Waltershausen quellenm\"a\ss{}ig zwar nur indirekt, aber durch
Zeitzeugenbericht konkret belegte Geschichte der Gau\ss{}schen
Euklidizit\"atsprüfung des physikalischen Raumes auf Grundlage der
geod\"atischen Pr\"azionsmessungen wurde 1972 durch A. Millers Artikel ohne Diskussion 
 des wichtigsten Zeugenberichtes und mit zweifelhaften  Argumenten  zu
einem ``Mythos'' umkonstruiert, erg\"anzt durch eine offensichtlich falsche
Hypothese über die Genese des ``Mythos''. Es könnte interessant sein, der
Frage nachzugehen,  welcher Konjunktion von Bedingungen  es geschuldet war,
dass dieser Vorsto\ss{} so erhebliche  Resonanz in der
Wissenschaftsgeschichte fand. Das  ist aber hier nicht unsere Absicht;
vielmehr möchte ich  im n\"achsten Abschnitt kurz vorstellen,  welche
methodologisch und historisch interessanten Fragestellungen trotz aller
M\"angel des ursprünglichen Anlasses aus der erst einmal  angeregten Debatte
um die Glaubwürdigkeit des Sartoriusschen Berichtes 
hervorgingen. Dabei geht es mir  hier nicht darum, die   Debatte in ihrer
ganzen Breite  zu verfolgen. Ich werde mich stattdessen auf eine Auswertung
des von E. Breitenberger im {\em Archive for the History of Exact Sciences}
publizierten Artikels \cite{Breitenberger:Gauss}  aus Sicht unserer
Fragestellung  konzentrieren, um zu lokalisieren, welche ernst zu nehmenden
Zweifel am Bericht des Sartorius  sich im Verlauf der
Debatte   herausbildeten. Der Kürze halber werde ich die These, Gau\ss{}
habe nie auf Grundlage  seiner geod\"atischen Pr\"asizisionsmessungen  des
$\triangle BHI$  eine Euklidizit\"atsüberprüfung des physikalischen Raumes
vorgenommen,  im folgenden als {\em  Millerschen Mythos} bezeichnen.

\subsubsection*{3. Weiterführende historisch-methodologische Zweifel }
In seinem  sachkundigen und sehr aufschlussreichen   Beitrag von
1984 stellte E. Breitenberger die neuere Kritik am ``Mythos'' der Gau\ss{}schen
Euklidizit\"atsüberprüfung in einen breiteren Rahmen. Er verwies auf \"altere
Bemühungen, des Astronomen Hugo von Seeliger und des Listing-Schülers Edmund
Hoppe, in \"ahnliche Richtung.\footnote{ \cite{Hoppe:Gauss}, \cite[614]{Kienle:Seeliger}.}
Allerdings  erw\"ahnte er dabei nicht, dass Seeliger durchaus eigene Gründe
dafür hatte,  Gau\ss{} von dem --- wie Seeliger meinte --- ``Makel''   befreit
sehen zu wollen,  sich   konkret (mit Suche nach jeweils bestmöglichen
quantitativen Fehlerschranken) um die  empirischen Gültigkeitsbedingungen
der euklidischen Hypothese gekümmert zu haben. Seeliger vertrat n\"amlich
selber in der Debatte des sp\"aten 19. Jahrhunderts um die nichteuklidische
Geometrie die Position, dass diese nicht einmal als  geometrische
Hintergrundstheorie   der Astronomie brauchbar sei. Er hatte bei seiner
Gau\ss{}-Interpretation also  durchaus Parteiliches  im
Sinne.\footnote{Breitenberger zitiert lediglich den Seeligerschen
``Makel''-Vorwurf, ohne ihn in den Kontext der Seeligerschen Auffassung zu stellen
\cite[274]{Breitenberger:Gauss}. Eine der Thesen, die K. Schwarzschild  in seiner Habilitation (München 1899, bei  Seeliger) zu verteidigen hatte,  beinhaltete  genau  diese  Auffassung einer Bedeutungslosigkeit der NEG für die Astronomie. Schwarzschild nahm dies zum 
Anlass für  eine höchst nuancierte  Diskussion der Frage in  \cite{Schwarzschild:1900}; vgl. dazu \cite{Schemmel:Schwarzschild}.}
Hoppe vertrat dagegen (unter Bezug auf einen Bericht Listings aus den 1870er
Jahren) die Auffassung, Gau\ss{} habe die empirischen Gültigkeitsschranken der
euklidischen Geometrie nicht aus geod\"atischen sondern aus astronomischen
Messungen zu bestimmen versucht. Ich werde im letzten Abschnitt darauf
zurückkommen.

Breitenberger steuert in seinem Artikel einen wendungsreichen Kurs durch die
widersprüchliche Literaturlage zu diesem Thema. Er stützt die Kritiker von
Seeliger bis Miller durch  Argumente, die er als {\em interne Evidenz}
(``internal evidence'') dafür betrachtet, dass  Gau\ss{} eine empirische
Überprüfung der euklidischen Raumstruktur durch seine geod\"atischen Messungen
aus methodologischen Gründen gar nicht h\"atte durchführen können.
Andererseits kennt Herr Breitenberger die Zeitzeugenberichte über eben
solche Überlegungen von Gau\ss{} zu gut, um sie wie in der einfachen Variante
des Miller-Mythos zu ignorieren oder wie in sp\"ateren ``Verfeinerungen'',
Sartorius als unglaubwürdigen weil mathematisch ignoranten Zeugen zu
entwerten. Er kommt  zu der von ihm im einzelnen  gut belegten
Schlussfolgerung:
\begin{quote}
\ldots  it is safe to conclude that Gauss was sufficiently irked by the
axiom of parallels, to bring it up in conversations repeatedly, and in
different forms, sometimes  quoting BHI [das {\em gro\ss{}e} Dreieck, E.S.]
sometimes not.

Thus the myth of the BHI triangle as a deliberate test of Euclidean geometry
appears as a fanciful embroidery upon indubitable facts, encouraged possibly
by reports of remarks made  by Gauss in his inner circle.
\cite[289]{Breitenberger:Gauss}
\end{quote}

Breitenberger kommt  zu der überraschenden Schlussfolgerung,  Gau\ss{} habe im
inneren Kreis über empirische Kontrollüberlegungen  der Gültigkeit der
euklidischen  Geometrie gesprochen, die er  gar nicht  habe durchführen
können. Diese Einsch\"atzung verbindet er jedoch  mit einer  {\em
Verteidigung}  des Sartorius-Berichtes, sogar einschlie\ss{}lich der dort
gemachten  quantitativen Angaben:
\begin{quote}
It also stands to reason that  Sartorius could hardly have misunderstood the
sense of Gauss' remarks, for although he was no mathematician, he was
thoroughly familiar with surveying and had produced detailed maps of the
Etna volcano, among others. What he recalls [die Genauigkeitsschranke der
Winkelsumme von ``zwei Zehntheilen einer Sekunde''  im gro\ss{}en Dreieck, E.S.]
is indeed quite consistent with the evidence above.
\cite[289]{Breitenberger:Gauss}
\end{quote}
In der ``evidence above'' berechnet Breitenberger auf Grundlage der
Gau\ss{}schen Messdaten die Abweichung der auf das
ebene Dreieck umgerechneten Winkelsumme des $\triangle BHI $ zu $0,642''$
\cite[284f.]{Breitenberger:Gauss}, eine Möglichkeit auf die  schon  \cite{Waerden:Myth} hingewiesen hatte.  Die von Gau\ss{} angebenenen Messdaten des gro\ss{}en
Dreiecks führen n\"amlich auf eine ``gemessene'' Winkelsumme von $180^{\circ}
\, 0' \; 14,211''$ auf der Sph\"are. Rückrechnung auf ein seitengleiches
ebenes Dreieck durch ``Elimination des sph\"arischen Exzesses''  führt auf den angegebenen Wert. Breitenberger kommentiert
dies wie folgt:
\begin{quote}
The spherical excess of BHI is stated by Gauss himself to have been
$14.85348''$; hence the closure error was $-0.642''$. He might well call this
fine result {\em vortrefflich} and it astonishes us a little that he did not
make more of it. \cite[285]{Breitenberger:Gauss}
\end{quote}
Der berechnete Wert ($0,6''$) stimmt grö\ss{}enordnungsm\"a\ss{}ig mit dem von Sartorius angebenen Wert
überein; keiner der hier erw\"ahnten Autoren sieht in der leichten Abweichung
``zwei'' statt sechs ``Zehnteile'' eine die Quelle entwertende
Differenz, und auch Herr Breitenberger beurteilt den Vergleich als ``quite  consistent''. Gau\ss{} selber erw\"ahnte den Wert $0,6''$ am 28. 12. 1823 in seiner Korrespondenz mit  Olbers \cite[II, 266]{Gauss:Olbers}.\footnote{Es w\"are ja durchaus naheliegend,  Sartorius bei der Memorierung
des von Gau\ss{} angegebenen Pr\"azisionsma\ss{}es eine Vertauschung des
Abweichungswertes  der Winkelsumme (0,6'') mit dem auf Einzelwinkel
umgerechneten  ($0,2''$) zuzutrauen. Aus Gau\ss{}' Sicht hatte die Umrechnung des Fehlers auf die Einzelwinkel durchaus Sinn; siehe  dazu Anmerkung  (34). Den Hinweis auf den Brief vom 28.12. 1823 an Olbers verdanke ich  E. Breitenberger. }

Dennoch sieht Herr Breitenberger  darin {\em keinen} Beleg für eine zutreffende
Darstellung des {\em gedanklichen Inhalts} der Gau\ss{}schen Bemerkung durch
Sartorius.
Er interpretiert n\"amlich die Gau\ss{}schen Pr\"azisionsdaten lediglich
unter dem Gesichtspunkt eines {\em nur pragmatisch} verstandenen {\em
Schlie\ss{}ungsfehlers} in den Messdaten eines geod\"atischen Dreiecks, der im
Rahmen der euklidischen (Raum-) Geometrie berechnet wird. Zur Unterstützung
dieser Interpretation verweist er auf die übliche Praxis der Geod\"aten (vor,
nach und bei Gau\ss{}), in der eine Qualit\"atssicherung der Messdaten durch
Winkelsummen in den ausgemessenen Dreiecken durch einen Vergleich mit der
erwarteten Summe ($180^{\circ} = \pi $ im ebenen Fall, $\pi  + F/ R ^2$ bei
Dreieckssfl\"ache $F$ und Erdradius $R$ im sph\"arischen Fall) üblich war.  Im genannten Artikel wird sogar darüberhinaus grunds\"atzlich  {\em bestritten}, dass es
für Gau\ss{} {\em überhaupt möglich gewesen sein könnte,} die Daten  {\em auch} als
ein  Genauigkeitsma\ss{} für die empirische Gültigkeit der euklidischen
Geometrie anzusehen. In diesem Sinne folgt Herr Breitenberger dort  einer durch
mathematisch-methodologische Überlegungen verfeinerten Variante des
Miller-Mythos und nimmt  die Überraschung darüber in Kauf, dass Gau\ss{} aus
seinen  Pr\"azisionsdaten angeblich,  {\em gegen} den sonst  als glaubwürdig angesehenen Bericht des Gau\ss{}-Freundes
Sartorius,   ``nicht mehr gemacht'' habe. Immerhin wurde die Pr\"azision der
Gau\ss{}schen Messungen bis zur Einführung von Lasermessmethoden  in den 1960er
Jahren nicht merklich übertroffen.\footnote{Zum Vergleich der Gau\ss{}schen
Messungen mit der sp\"ater  erreichten Pr\"azision von Theodolitmessungen bis ca
1960, siehe die Verweise in \cite[Anm. 35, 36, 42, 50,
52]{Breitenberger:Gauss}.
Aufschlussreich ist dabei insbesondere die über Anm. 52 zitierte Bemerkung
aus einem   Geod\"asielehrbuch   der frühen 1970er:  ``\ldots angles measured
by \ldots theodolites beween 1800 and 1950 are not necessarily inferior to
those measured by modern instruments \ldots the long time required to
complete observations \ldots probably resulted in better elimination of
refraction errors than modern observers have time for '' (G. Bomford {\em
Geodesy.} Oxford 3rd ed. 1971, 18, zitiert nach
\cite[280]{Breitenberger:Gauss}).}

Herr Breitenberger stützt sich bei seiner Argumentation auf die ``interne
Evidenz'', dass Gau\ss{} in allen seinen geod\"atischen Arbeiten, selbst bei
seinen publizierten differentialgeometrischen Grundlagenuntersuchungen
(Differentialgeometrie der Fl\"achen als ``Grundlage der Grundlage'' der
Geod\"asie) die euklidische Geometrie des einbettenden Raumes stets
voraussetzte. Das warf natürlich für jeden Versuch einer empirischen
Kontrollüberlegung oder gar für eine Bestimmung von Genauigkeitsschranken
für die Gültigkeit der euklidischen Geometrie grunds\"atzliche Probleme auf.
In der Tat schuf ja erst Riemann  in seinem Habilitationsvortrag 1854 den
begrifflichen Rahmen für eine verallgemeinerte Geometrie, in dem  solcherart
Fragen systematisch zu kl\"aren
waren.\footnote{\cite[285]{Breitenberger:Gauss}.}
Nun ist aber wohlbekannt, dass Gau\ss{} dies nicht daran hinderte,
Grundbeziehungen der nichteuklidischen Geometrie zu erforschen, auch {\em
ohne} vorg\"angig den gesamten theoretischen Rahmen gekl\"art zu haben.

Zentrales methodisches Element war für Gau\ss{}  die Beziehung zwischen
charakteristischer Konstante $C$,  Krümmung (extrinsisch definiert,
intrinsisch bestimmbar) $\kappa$ einer Gau\ss{}schen
Fl\"ache und Winkelsumme im Dreieck im Fall  konstanter Krümmung
\beq \alpha +\beta +\gamma   - 2\pi   =  \kappa   F ,\eeq
sowie $ \kappa = - C $ im Fall der nichteuklidischen Geometrie. Gl. (1) war (und ist)  eine
Spezialisierung des Gau\ss{}schen  {\em Theorema elegantissimum}, d.h. des von
Gau\ss{} stammenden Anteils des sp\"ater nach  Gau\ss{}-Bonnet benannten  Satzes für den Fall
konstanter Krümmung und geod\"atischer Dreiecke. Es ist nicht anzunehmen, dass
Gau\ss{} etwa gezögert haben könnte, diese Beziehungen auch für die dreidimensionale
nichteuklidische Geometrie {\em heuristisch} als gültig anzusehen, obgleich eine
systematische Begründung von ihm nicht ausgearbeitet wurde und eine Ausarbeitung ohne die Riemannsche Verallgemeinerung seines differentialgeometrischen Ansatzes  die bekannten Schwierigkeiten zu überwinden gehabt h\"atte.  Seine Überlegungen zu den Grundlagen der Geometrie hatten ihn ja auch dort auf dieselbe Winkelsummenbeziehung geführt (nur eben mit der berühmten ``Constante'' C). Die 
eingangs zitierte Bemerkung  im Brief an Taurinus über
die Möglichkeit der Überprüfung der Geometrie durch ``Beobachtungen auf der
Erde oder am Himmel'' (wie auch andere \"ahnlich gelagerte) sprechen in dieser
Hinsicht eine klare Sprache.

Breitenberger tr\"agt m.E. dieser heuristischen Konstellation nicht
ausreichend Rechnung, wenn er die richtige (wohlbekannte) Beobachtung, dass
erst Riemann in ausgearbeiteter Form ``the notion of a curved
threedimensional space'' eingeführt h\"atte (p. 285) zu folgender
apodiktischer Behauptung zuspitzt:
\begin{quote}
Gauss was not alluding to a test of (three-dimensional) geometry; as
stressed above, he did not and he indeed could not, envisage one. He wrote
in ignorance of any tangible application of hyperbolic trigonometry to some
surface or object. Not only did he not possess a model in the modern sense
from which he might have deduced self-consistency;  he did not even possess
any realization which would have made the subject accessible to geometric
intuition. \cite[287]{Breitenberger:Gauss}
\end{quote}
Natürlich kannte Gau\ss{} keine ``hyperbolische Geometrie'' usw.; dennoch behaupte ich, dass Gau\ss{} mit den von ihm ausgearbeiteten und
favorisierten Konzepten der geod\"atischen Linie, einer intrinsisch durch
Metrik und Krümmung bestimmten Fl\"ache und den Winkelsummensatz im Spezialfall
konstanter Kümmung, entscheidende methodisch-begriffliche  Bausteine  zu
Verfügung standen, die es ihm möglich machten, Pr\"azisionsmessungen ``auf der
Erde oder am Himmel''   {\em heuristisch gut begründet} bezüglich  ihrer
Aussagef\"ahigkeit über die Natur des physikalischen Raumes  auszuwerten. Dies schloss die  Möglichkeit ein, pr\"azise Genauigkeitsschranken für die Gültigkeit der
euklidischen Geometrie zu bestimmen und anzugeben. {\em Es bedeutet dabei keine starke Zuspitzung, eine
solche Genauigkeitsschranke als  informelle Fassung der Absch\"atzung einer
mit den Messergebnissen vertr\"aglichen oberen  Schranke für den Betrag der
Raumkrümmung zu interpretieren.} In diesem Sinne halte ich die oben zitierte
Bemerkung, ``he did not even possess any realization which would have made
the subject accessible to geometric intuition'' für falsch. Natürlich h\"angt
bei der Beurteilung einer solchen Aussage alles von der Interpretation
``geometric intution'' ab; sicher  lie\ss{}en sich auch Pr\"azisierungen angeben,
die Gau\ss{} nachweislich nicht zur Verfügung standen, historisch aber nicht
erhellend sind. Eine Diskussion auf dieser Ebene scheint mir 
fruchtlos.

Stattdessen möchte ich  den  von Herrn Breitenberger  recht apodiktisch
ausgesprochenen Zweifel in eine  Serie von Problemen auflösen,  die zu
kl\"aren sind, damit die  Frage  entscheidbar wird, ob die  geod\"atischen Pr\"azisionsmessungen des $\triangle BHI $   im Rahmen der Gau\ss{}schen Methodologie heuristisch (und im Rahmen der Riemannschen Erweiterung exakt),  als
Genauigkeitsschranken für die empirische Gültigkeit der euklidischen
Geometrie ausgewertet werden konnten. Jeremy Gray hat mich auf die meisten
der hier zu kl\"arenden Fragen durch beharrliches Nachfragen aufmerksam
gemacht. Den Diskussionen mit ihm verdankt dieser Beitrag viel; genauer
gesagt w\"are er ohne diese Diskussionen nicht geschrieben
worden.\footnote{Ein Dank geht auch an P. Strantzalos und M. Lambrou, denen
J. Gray und ich die  Gelegenheit verdanken, unsere damals noch
divergierenden Auffassungen  auf  einem   Workshop über $\Sigma \upsilon \nu
\acute{\epsilon} \delta \rho \iota o \; I \sigma \tau o \rho \dot{\iota}
\alpha \sigma  \;\kappa \alpha \iota  \; \Delta \iota \delta \alpha \kappa
\tau \iota \kappa \acute{\eta} \sigma \; \tau \omega \nu \; M \alpha
\theta \eta \mu \alpha \tau \iota \kappa \acute{o} \nu  $ im April 2002 an
der Universit\"at Kreta (Heraklion)  in anregender und  freundschaftlicher
Atmosph\"are ausführlich zu diskutieren.}

Der Übersichtlichkeit halber werde ich  in Teil 4 die grunds\"atzlichen
Probleme diskutieren, die bei dem Versuch zu berücksichtigen sind, die
Raumkrümmung durch Messungen ``auf der Erde'' zu bestimmen, ohne mich den
Restriktionen der Gau\ss{} zur Verfügung stehenden, fertig ausgearbeiteten
mathematischen  Theorien zu unterwerfen. Im nachfolgenden Abschnitt (Teil 5)
folgt dann eine Diskussion aus  der Gau\ss{}schen Sicht der 1820er Jahre. Des
Kontextes wegen bitte ich die Lesenden,  bezüglich der empirischen  Geltung
der euklidischen Geometrie lediglich das abgesicherte Wissen des frühen 19.
Jahrhunderts in Rechnung zu stellen. Astronomische Parallaxmessungen verschwanden
zu dieser Zeit noch im ``Schmutz'' der Messeffekte. Die grö\ss{}ten gut
vermessenen geod\"atischen Dreiecke hatten Seitenl\"angen in der Grö\ss{}enordnung
mehrer 10 km, die Standardabweichungen der Richtungsmessungen der von Gau\ss{}
als Bezug dienenden niederl\"andischen Netzmessung des Baron C. von Krayenhoff
lag bei $2,7''$,\footnote{Nach Kalkulationen von Gau\ss{} auf Basis der ihm
bekannten  Daten; siehe \cite[281]{Breitenberger:Gauss}. Dort ist
f\"alschlich vom d\"anischen Gradnetz die Rede; das stand allerdings unter
Leiting von  H.C. Schumacher.}
Bei einem grö\ss{}eren Dreieck dieses Netzes lag der  Fl\"acheninhalt  um die  300
km$^2$  und damit der  sph\"arische Winkelexzess bei 1,5'',  also {\em  unter}
dem ``mittleren Fehler'' (Standardabweichung) einer Richtungsmessung. Schon
der sph\"arische Exzess eines einzelnen (direkt vermessenen) Dreiecks
verschwand tief im Schwankungsbereich  der Messfehler.

 Aus Sicht unserer Geschichte l\"asst sich die damit erreichte Pr\"azision
pointiert  so umschreiben: W\"are jemand vor den
Gau\ss{}schen Pr\"azisionsmessungen auf die Idee gekommen, eine vermutete
Raumkrümmung abzusch\"atzen, so h\"atte dies (zumindest stillschweigend) die
Annahme vorausgesetzt, dass die Raumkrümmung  die Krümmung der
Erdoberfl\"ache {\em erheblich überschreitet} (!); sonst w\"are bei der bis
dahin erreichten  Messgenauigkeit überhaupt kein Effekt für
Winkelsummenmessungen  beobachtbar geworden. Diese, wie ich meine, nur
scheinbar  paradox  klingende Umformulierung der vor Gau\ss{} erreichten
Messgenauigkeit sollte deutlich machen, worum es ging, wenn man --- wie Gau\ss{}
--- zu Beginn des 19. Jahrhunderts einen konsequenten Empirismus in Sachen
metrischer Raumstruktur zu verfolgen beabsichtigte. Dabei war  ja   {\em
nicht auszuschlie\ss{}en},\footnote{Dies ist offensichtlich etwas anders, als
(positiv) den entsprechenden Sachverhalt zu vermuten.}
dass die vermutete Konstante ``in einigem Verh\"altnisse zu (\ldots) Grö\ss{}en
(\ldots) im Bereich unserer Messungen auf der Erde''  stehen könnte. Wie
sollte man dann aber ---  wenn man die a-priorische Gültigkeit der euklidischen Geometrie nicht akzeptierte --- anders als empirisch sicherstellen, dass  bei einer Erhöhung der Messpr\"azision
zusammen mit den Effekten der Erdkrümmung keine Störungen auftreten können,
die möglicherweise von der Raumkrümmung herrührten?

\subsubsection*{4.  Bestimmung von Schranken für die  Raumkrümmung aus
erdgestützten  Messungen }
Natürlich waren bei einem  Projekt der Absch\"atzung der Raumkrümmung durch
terrestrische Messungen einige wohlbekannte Sachverhalte in Rechnung zu
stellen:
\begin{Sachverhalt}
Zu jedem Wert  $s$ mit  $\pi < s  < 2\pi $ gibt es ein sph\"arisches Dreieck zugehöriger Winkelsumme  $s  = \alpha + \beta +\gamma $ .
\end{Sachverhalt}

Diese bekannte Tatsache allein schlie\ss{}t schon aus, die Methode der
Krümmungsbestimmung durch  Winkelsummenmessung in einem ebenen Dreieck
direkt (das hei\ss{}t, ohne Einbeziehung weiterer metrischer Daten der
Erdoberfl\"ache) auf geod\"atische Dreiecke der Erdoberfl\"ache zu übertragen;
letztere können  ja nach
\cite[ \S 28]{Gauss:1828} in ausreichender Genauigkeit als sph\"arische
Dreiecke betrachtet werden.  

Aber auch die Einbeziehung intrinsischer
metrischer Daten der Fl\"ache l\"asst nicht auf Anhieb einen Zugang zur
Absch\"atzung der Raumkrümmung zu. Dies liegt an:
\begin{Sachverhalt}
Lokale Messungen (der Metrik) in einer 2-dimensionalen Untermannigfaltigkeit
einer 3-dimensionalen Riemannschen Mannigfaltigkeit lassen keinen
Rückschluss auf die Krümmung des einbettenden Raumes zu.
\end{Sachverhalt}

Man denke etwa nur an die Horosph\"are des nichteuklidischen Raumes, die lokal
euklidische Geometrie aufweist, oder an eine Gau\ss{}sche Fl\"ache im euklidischen
Raum,  die in ihrer intrinsischen Geometrie nicht von einer isometrisch in
einen nichteuklidischen Raum oder von  einer in eine allgemeinere
Riemannsche Mannigfaltigkeit eingebetteten Fl\"ache unterscheidbar ist usw..
Darauf spielt Breitenberger verschiedentlich an, wenn er betont, dass die
Gau\ss{}schen Berechnungen stets die euklidische
Geometrie voraussetzten, anstatt sie zu überprüfen. Eine Überprüfung hat ja
nur dann Sinn, wenn ihr Ausgang  grunds\"atzlich auch negativ sein
kann.

Dieser Sachverhalt wird, in unterschiedlichen Formulierungen vorgebracht,
von verschiedenen Autoren   als eine schlüssige {\em Widerlegung der
Möglichkeit} überhaupt angesehen, eine Überprüfung der euklidischen
Raumstruktur durch terrestrische Messungen ausführen zu können. Wie sich
gleich zeigen wird, ist dies jedoch {\em falsch}. W\"are diese Ansicht
richtig, so beruhten  Gau\ss{}ens Bemerkungen aus den frühen 1820er Jahren über
die Möglichkeit, eine Abweichung des physikalischen Raumes von der
euklidischen Geometrie eventuell durch Messungen ``auf der Erde''
feststellen zu können, tats\"achlich auf einem Irrtum.\footnote{Man denke etwa
an Seeligers ``Makel'' auf dem Bild des wissenschaftlichen ``Heroen'' Gau\ss{}.}
 Man schlie\ss{}e jedoch nicht zu schnell, sondern vergesse nicht: 
\begin{Sachverhalt}
Kann man erdgestützt die Winkelsumme eines gro\ss{}en ebenen extrinsischen (d.h.
nicht in der Erdoberfl\"acher liegenden) Dreiecks  messen oder diese aus den
direkt erfolgenden Messungen rückrechnen, so erh\"alt man aus der Gau\ss{}schen
Winkelsummenbeziehung (1) eine direkte oder indirekte Information über die Raumkrümmung.
\end{Sachverhalt}

Das hier als ``extrinsisch'' bezeichnete Dreieck kann  durch
Lichtstrahlenverbindungen zwischen Bergspitzen (etwa $\hat{B}, \hat{H},
\hat{I}$)) realisiert  werden, die über die Oberfl\"ache $\cal{F}$ der
mathematischen Erdgestalt  (Sph\"are, Ellipsoid oder Geoid) hinausragen und
direkte Sichtverbindung untereinander zulassen. Das Dreieck kann  und wird
im allgemeinen schief, also nicht parallel  zu den drei Tangentialebenen an
$\cal{F}$ in den orthogonal unter den Bergspitzen liegenden Punkten $B, H, I
\in \cal{F}$ liegen. Die Orthogonalit\"at kann als durch Schwerkraftwirkung
definert angesehen werden (Vertr\"aglichkeit der  Lichtstrahlengeometrie mit
der Orthogonalit\"atsbeziehung bezüglich Schwerefeld
unterstellt).\footnote{Die Orthogonalit\"atsbeziehung wird also genau genommen
zum Geoid betrachtet. Dies kann aber für die folgenden Überlegungen  in
ausreichender Approximation mit dem Referenzellipsoid zusammenfallend
angesehen werden.}
 Weiter wird hier auf  Lichtstrahlen als bestmögliche empirische
Realisierungen für Geod\"atische Bezug genommen. Dies schlie\ss{}t Korrekturen bei
der Messauswertung, die Beugung durch Lufteinflüsse auf die
Lichtausbreitung als Fehlerquelle betrachten und zu kompensieren versuchen,
ausdrücklich mit ein und entspricht  einer treuen Interpretation der
Gau\ss{}schen Praxis (und der sp\"aterer Akteure).\footnote{Vergleiche dazu auch
\cite[282]{Breitenberger:Gauss}.}  Genau genommen geht und ging es also  um
eine so weit wie möglich {\em von Atmosph\"areneinflüssen bereinigte Geometrie der
Lichtstrahlen}.
Die Höhen von  $\hat{B}, \hat{H}, \hat{I}$ über dem von Gau\ss{} verwendeten
Referenzellipsoid betrugen  jeweils   1156 m (Brocken), 508 m (Hohehagen) und 916 m
(Inselsberg).\footnote{Werte angegeben nach
 \cite[279 ,283]{Breitenberger:Gauss}.}
\newline{}

\unitlength1cm
\begin{figure}[t]
\begin{picture}(10,7)
 \includegraphics[height=7cm]{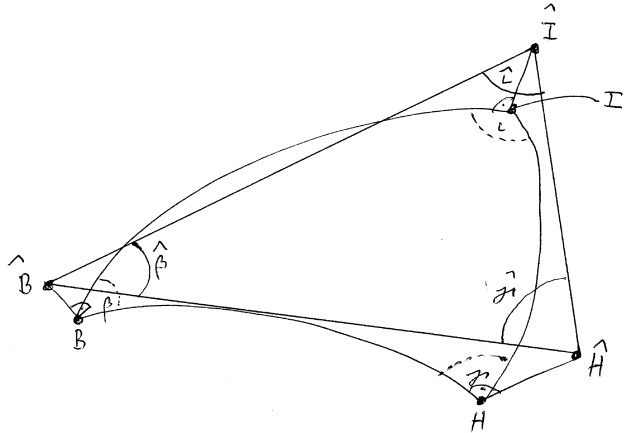} 
\end{picture} 
\caption{Geod\"atisches Dreieck und Lichtstrahlendreieck}
\end{figure}

Nun wirkt sich eine mögliche Abweichung der extrinsischen
 Lichtstrahlen-Geometrie vom euklidischen Fall auf die Winkelsumme im $\triangle \hat{B}
\hat{H} \hat{I} $ aus und damit auch auf die bei $B, H, I$ in die jeweiligen 
Tangentialebenen an $\cal{F}$ projizierten Winkel und deren Summe. Gemessen
werden aber in der geod\"atischen Praxis nicht die in der Ebene durch
$\hat{B}, \hat{H}, \hat{I}$ liegenden Winkel $\hat{\beta }, \hat{\gamma },
\hat{\iota }$ des  Dreiecks $\triangle\hat{B} \hat{H} \hat{I} $, sondern die
orthogonal in die Tangentialebenen an $\cal{F}$ projizierten Winkel $\beta
,\gamma ,\iota $. Dies ist ein Ergebnis des Messverfahrens selber und
erfordert keine zus\"atzlichen Rechnungen: Nach einer ``Horizontierung des
Theodolits'' (das hei\ss{}t Einrichtung des Instrumentes in die Horizontalebene)  wird die
Winkeldifferenz (``Azimut'') der Peilrichtung (genau genommen der
Vertikalebene durch den Peilstrahl) zur Südrichtung abgelesen. Die Winkel
des geod\"atischen Dreiecks $BHI$ ergeben sich als Differenz der Azimute der
jeweiligen Richtungen $BH, BI$ usw..\footnote{Diese Berechnung der Winkel
aus den Gau\ss{}schen Azimutdaten beschreibt auch Breitenberger auf seiner p.
284.  Unglückseligerweise gibt er  an anderer Stelle eine irrtümliche
Beschreibung der Arbeitsg\"ange, die Gau\ss{}  bei der Auswertung der Messdaten
angeblich ausführen musste. Dazu z\"ahlt er vier Arbeitsschritte auf, von
denen er die ersten beiden wie folgt beschreibt: ``The data first had to be
reduced to uniform station height; for this he  used formulae of his own
(\ldots). Then the lines of sight, which form the edges of a polyhedron
inscribed into a sphere, were projected onto the sphere from its center,
becoming arcs of great circles, and the sum of the angles in the resulting
spherical triangles had to be increased over $180^{\circ}$ by the spherical
excess which for each triangle would be calculated form its size....'' Es
folgen zutreffende Beschreibungen des Fehlerausgleichs und der Berechnung
von Koordinaten
\cite[280]{Breitenberger:Gauss}. Da für die gesamte weitere Diskussion entscheidend ist, welche Daten direkt gemessen werden und welche durch theoretische Rückrechnung zu erhalten sind, ist hier deutlich darauf hinzuweisen, dass  bei den beiden ersten
 hier aufgeführten Schritten  theoretische Beschreibung
und Auswertung in irreführender Weise  miteinander vermengt werden. Tats\"achlich erfolgt ja die Peilung l\"angs der Kanten
eines Polyeders, dem auf die hier beschriebene Weise ein der Erdgestalt
einbeschriebenes Polyeder zugeordnet werden kann (aber nicht muss). Die
Winkelbestimmung wird bei dieser Beschreibung aber g\"anzlich auf
den Kopf gestellt. Durch Differenzbildungen der direkt gemessenen
Azimutwerte werden --- wie hier im Haupttext beschrieben und von E. 
Breitenberger im Fall der Gau\ss{}schen Werte für das Dreieck BHI sogar
vorgeführt (!) --- die Winkel eines aus dem Polyeder-Dreieck auf die Fl\"ache
projizierten geod\"atischen Dreiecks  bestimmt. Der sph\"arische Exzess wird ---
wie von Gau\ss{} vorexerziert --- umgekehrt eingesetzt, um die Winkel eines
seitengleichen ebenen Dreiecks zu bestimmen (von Gau\ss{} verallgemeinerter Satz
von Legendre)! Schon der erste in Breitenbergers Text  angegebene Schritt scheint
auf einem Missverst\"andnis zu beruhen: Die Projektion der Winkel in die
Tangentialebene ist wie oben beschrieben gewisserma\ss{}en ``in den Messprozess
eingebaut'', eine Reduktion auf gleiche Stationshöhe ist für die Bestimmung der Azimute und damit der geod\"atischen Winkel nicht erforderlich.   Ich danke Herrn Heindl, numerischer Mathematiker und gelernter Geod\"at, für die ihm leicht fallenden, in diesem Kontext \"au\ss{}erst hilfreichen Erl\"auterungen über 
Theodolitmessungen.}

Es ist nun abzusch\"atzen, wie sich eine kleine Abweichung der extrinsischen
Raumgeometrie von der euklidischen Metrik, etwa mit (konstanter Schnitt-)
Krümmung $\hat{\kappa }$ ( $|\hat{\kappa }| << 1$,  für Gau\ss{} $\hat{\kappa
}<0$), auf die geod\"atischen Winkelmessungen auswirken würde. Wir gehen zur
Vereinfachung davon aus,
 dass durch unabh\"angige  {\em intrinsische} Messungen die metrischen Daten
der Erdoberfl\"ache $\cal{F}$ bekannt seien (mittlere Krümmung,
Hauptkrümmungen des Ellipsoids), sowie die Koordinaten der Dreieckspunkte
$B, H, I$ und damit der Fl\"acheninhalt $F(\triangle )$  des durch Projektion der
Lichtstrahlenverbindungen gewonnenen geod\"atischen Dreicks $\triangle
BHI$.\footnote{Zur Diskussion, inwieweit diese Vereinfachung durch die bei
Gau\ss{} vorliegende Messpraxis abgedeckt erscheint, siehe unten.}
Die Winkel $\hat{\beta }$ usw. weichen von einem bei euklidischer
Raumgeometrie angenommenen Wert $\hat{\beta _ {0}}$  nur um einen kleinen
Betrag $\Delta \hat{\beta }$ etc.  ab, die Winkelsummenabweichung sei
$\epsilon $:
\[ \hat{\beta } = \hat{\beta _ {0}} + \Delta \hat{\beta }, \; \;
\hat{\gamma  } = \hat{\gamma  _ {0}} + \Delta \hat{\gamma  } , \; \;
\hat{\iota  } = \hat{\iota  _ {0}} + \Delta \hat{\iota  } \]
\[ \hat{\beta _ {0}} + \hat{\gamma  _ {0}} + \hat{\iota  _ {0}} = \pi , \;
\;
 \Delta \hat{\beta } + \Delta \hat{\gamma  } + \Delta \hat{\iota  } =:
\epsilon  \]

Eine elementargeometrische Überlegung und Reihenentwicklung zeigt, dass die
 unterstellten nichteuklidischen Winkeldeformationen $\Delta \hat{\beta }$
usw. sich  bei der Orthogonalprojektion in erster N\"aherung mit einem nahe
bei 1 liegendem Faktor als Änderungen der Winkel gemessenen geod\"atischen
Winkel $\beta , \gamma , \iota $ auswirken (siehe Anhang, Gl. (5)). Es ergibt
sich also:
\begin{Sachverhalt}[a]
Eine kleine Abweichung der Winkelsumme im Dreieck $\hat{\triangle } :=
\triangle\hat{B} \hat{H} \hat{I} $ um $\epsilon = \hat{\kappa}\cdot
F(\hat{\triangle} ) $ vom euklidischen Wert $\pi $ wirkt sich in guter
erster N\"aherung als gleichgro\ss{}e Abweichung der Summe der auf der Sph\"are gemessenen Winkel $\beta , \gamma , \iota $  vom erwarteten sph\"arischen Wert aus:
\beq
\beta +\gamma +\iota  - \pi \approx \frac{1}{R^2} \cdot F(\triangle ) + \epsilon 
\eeq
($R$ mittlerer Erdradius, $ F(\triangle )$  Fl\"acheninhalt des Dreiecks $\triangle $).
\end{Sachverhalt}
 In umgekehrter Richtung folgt:
\setcounter{Sachverhalt}{3}
\begin{Sachverhalt}[b]
 Eine angenommene extrinsische Raumkrümmung $\hat{\kappa }$ wirkt sich bei
Vermessung ``gro\ss{}er'' geod\"atischer Dreiecke $\triangle $ als systematischer
Anteil $\epsilon '$ des ``Schliessungsfehlers''  der auf den euklidischen
Fall umgerechneten Winkelsumme aus. Dieser wird beim Fl\"acheninhalt $F(\triangle )$ des
Dreiecks in erster N\"aherung   abgesch\"atzt durch:\footnote{Hierbei wurde als weitere N\"aherung
$F(\hat{\triangle}) \approx  F(\triangle)$ verwendet.}
\beq \epsilon ' \approx \hat{\kappa} \cdot F(\triangle) \eeq
\end{Sachverhalt}

Wie schon Ende des letzten Abschnittes erw\"ahnt, lag der Schlie\ss{}ungsfehler
bei Messungen geod\"atischer Dreiecke bis etwa 1820 {\em über} dem sph\"arischen
Exzess. Allein schon deswegen w\"are es unmöglich gewesen, auf der vorliegenden
empririschen Grundlage eine Absch\"atzung der Raumkrümmung wie in (3) auch nur
zu erw\"agen.  Hinzu trat, dass die oben gemachte Voraussetzung (Gewinnung der
metrischen Erddaten durch intrinsische Messungen, Bestimmung der Koordinaten
der Dreieckspunkte durch unabh\"angige Messungen) für die vor-Gau\ss{}schen
Messungen  kaum  als erfüllt angesehen werden konnten, da das zur Prüfung
verwendete Dreieck Teil des Bezugsnetzes  und von derselben Grö\ss{}e wie die
Netzdreiecke war.

Bei  der Gau\ss{}schen Messung verhielt sich das anders: Das ``gro\ss{}e'' Dreieck
BHI hatte Seitenl\"angen (69 km (BH), 85 km (HI), 107 km (BI)) und war damit
den Seitenl\"angen nach um etwa den Faktor 5, der Fl\"ache nach um eine ganze
Grö\ss{}enordnung grö\ss{}er als die üblichen Netzdreiecke geod\"atischer Messungen,
einschlie\ss{}lich jener, die zur Bestimmung der Erdgestalt gedient hatten (oder
Gau\ss{} zur Datenkorrektur dienten).
Ging man davon aus, dass bei den {\em kleineren Netzdreiecken} auch bei den von
Gau\ss{} verbesserten Messmethoden ein eventueller Einfluss der Raumkrümmung
nicht nachzuweisen war, so konnte man die daraus abgeleiteten metrischen
Werte mit gutem Grund als zur {\em intrinsischen Geometrie} der Erdoberfl\"ache im
Sinne der Gau\ss{}schen Fl\"achentheorie gehörend ansehen. Die Peillinien des
Netzes, die im einbettenden Raum verliefen,  lagen ja offenbar in so kleinen
Stücken des Raumes, dass letztere ausreichend gut euklidisch darstellbar waren.  In der Sprache der Riemannschen Mannigfaltigkeiten hie\ss{}e das: Jedes
einzelne Lichtstrahlendreieck war  ausreichend genau (d.h. im Rahmen der
Messgenauigkeit)  in einem Tangentialraum an die Mannigfaltigkeit darstellbar,
  das darunter liegende geod\"atische Dreieck entsprechend auf einer euklidischen Sph\"are zu betrachten.

  Darüberhinaus hatte Gau\ss{} Wert darauf gelegt, das gro\ss{}e Dreieck individuell
mit höchster Genauigkeit zu vermessen und die Bestimmung seiner Grunddaten
nicht in den Netzausgleich des nach Hamburg verlaufenden eigenen Netzes,
des südlich anschlie\ss{}enden seines kurhessischen Kollegen
Gerling oder  des im Westen liegenden niederl\"andischen Netzes des Barons von Krayenhoff einzubeziehen.\footnote{\cite{Gerardy:Gauss}.}
 Methodisch hatte Gau\ss{} also darauf geachtet, dass das Dreieck au\ss{}erhalb der
Referenznetze lag und in diesem Sinne diesem ``methodisch  extrinsisch''
war. Insofern ist die Überlegung von {\em Sachverhalt 4a} und {\em 4b} auf
das Gau\ss{}sche $\triangle BHI$ (bzw.  $\triangle\hat{B} \hat{H} \hat{I} $)
begründet anwendbar.

Mit Fl\"acheninhalt $F(\triangle ) \approx  2920 $ km$^2$ (aus Gau\ss{}' Wert
$14,86''$ für den sph\"arischen Winkelexzess und seinem Wert für den Erdradius $R \approx 6370$ km
berechnet) und Schlie\ss{}ungsfehler $\epsilon ' \approx 0,6''$ gibt eine einfache Überschlagsrechnung\footnote{Etwa
$F/R^2 \approx 14''$, $F/r^2 < 0,6'' $, also $r/R > \sqrt{14/0,6} \approx 5
$.}
  für die Raumkrümmung $ |\hat{\kappa }| = 1/r^2$ die Schranke
\beq
r  > 5 \, R \approx 3 \cdot10^4 \;
\mbox{km}   \; \;
\mbox{beziehungsweise } \; \;\; 
| \hat{\kappa }|< 10^{-9} \; \mbox{km}^{-2} .
\eeq

 Es war also durchaus möglich,  eine obere Schranke für eine mit den
Gau\ss{}schen geod\"atischen Messungen vertr\"agliche Raumkrümmung zu bestimmen. Das
lag nicht nur an der erhöhten Genauigkeit um etwa eine Grö\ss{}enordnung,
sondern auch daran, dass Gau\ss{} seine Messungen so organisiert hatte, dass eine
Kontrollüberlegung dieser Art überhaupt  methodisch gerechtfertigt
erschien.  Sehr viel mehr war von  terrestrischen Messungen kaum   zu
erwarten.\footnote{Wie \cite[282]{Breitenberger:Gauss} es ausdrückt:
``[Gauss] had by 1823 attained geodetic standards which were essentially
unsurpassed one-and-a-half century later''.}
 Wenn der physikalische Raum gekrümmt sein sollte, so lag nun  die untere
Schranke für den Krümmungsradius jedenfalls merklich jenseits  des
Erdradius: $ r > 5 R$. Das erscheint aus astronomischer Sicht als l\"acherlich
gering; man berücksichtige aber, dass bis 1838 keine verl\"asslichen Werte für
Parallaxmessungen vorlagen und daher jeder Versuch einer quantitativen
Auswertung astronomischer  Messungen für Absch\"atzungen der Raumkrümmung
fast wörtlich zu lesen in der Luft  (genauer im ``Äther'') hing.

Diese Auswertung hat bis hierher theoretische Methoden und Begriffe
herangezogen, die Gau\ss{}  nicht ausformuliert  zur Verfügung standen. Es
bleibt also zu diskutieren, welche Aspekte der hier ``modernisiert''
vorgetragenen (allerdings schon ab 1854 durchführbaren) Analyse ihm auf
seine  Weise zug\"anglich gewesen sein konnten.

\subsubsection*{5. Was war für Gau\ss{} absch\"atzbar? }
So klar es ist, dass Gau\ss{} in den 1820er Jahren {\em in ausgearbeiteter Form}
keine 3-dimensionalen Riemannschen Mannigfaltigkeiten kannte (nicht einmal
konstanter Krümmung), so klar zeigt  seine  Korrespondenz, dass er seit etwa
1816 von der logisch-begrifflichen Möglichkeit einer konsequenten
nichteuklidischen Fassung der Geometrie  überzeugt
war.\footnote{\cite{Staeckel:Gauss}, siehe auch \cite[25ff.]{Reichardt:Gauss}} Das ist in der Wissenschaftsgeschichte
unumstritten. Auch ist bekannt, dass ihn  dabei insbesondere  interessierte,
wie die bekannten metrischen Beziehungen der euklidischen Geometrie, etwa
die des Winkelsummensatzes, modifiziert werden.
Einen Satz wie den oben zitierten an Taurinus  (``W\"are die Nicht-Euklidische
Geometrie die wahre, und jene Constante in einigem Verh\"altnisse zu solchen
Grö\ss{}en, die im Bereich unserer Messungen auf der Erde oder am Himmel liegen,
so lie\ss{}e sie sich a posteriori ausmitteln''), kann man daher   im Sinne des
obigen  {\em Sachverhalts 3}  interpretieren, ohne die Gau\ss{}sche Perspektive
auch nur einen Deut zu überschreiten. Kurz formuliert:

{\em Sachverhalt 3 war für Gau\ss{} klar.}

Ebenso bedurfte Gau\ss{} keineswegs einer Einordnung in die entfaltete
Theorie Riemannscher Mannigfaltigkeiten um zu sehen, dass die NEG im Kleinen
durch die euklidische Geometrie approximiert wird, insbesondere und umso
genauer je kleiner die  ``Constante''  $ C = |\kappa | << 1$ ist, egal ob im
ebenen oder im r\"aumlichen Fall.

Damit war aus seiner Perspektive völlig naheliegend,  die Vermessung der
Erddaten auch dann in guter Approximation als Sachverhalt der
intrinsischen Geometrie der Erdoberfl\"ache   anzusehen, wenn er für
Kontrollüberlegungen der empirischen Grundlagen die Geometrie des
physikalischen Raumes als nichteuklidisch, aber eben mit kleinem $C$, ansah.
Dies ist keine methodologische Widersprüchlichkeit, sondern legitime  und
gut gestützte Vorgehensweise empirischer mathematischer Naturwissenschaft.

Da im ``sehr'' Kleinen für die Raumgeometrie auch bei Annahme $C>0$ mit
gro\ss{}er Genauigkeit die euklidische Geometrie gilt, war und ist es 
selbstverst\"andlich, dass  die Projektion der Winkel des
Lichtstrahlendreiecks $\triangle\hat{B} \hat{H} \hat{I} $  euklidisch
approximiert werden kann.\footnote{Die entsprechende
Approximationsüberlegung ist natürlich  auch ohne Tangentialr\"aume an
Mannigfaltigkeiten  durchführbar. Sie  muss dann eben nur anders (aus
unserer Sicht etwas umst\"andlicher) formuliert werden Bei der Projektion in
den Eckpunkten des Dreiecks geht es um die Darstellung r\"aumlicher Gebiete in
der Grö\ss{}enordnung von jeweils 1 m (in alle drei Raumdimensionen). Es reicht
ja, die Projektion in eine Parallelebene zur Tangentialebene an das Geoid
durch den Messpunkt vorzunehmen, in dem der Theodolit steht.}  Dass dabei
eine kleine Winkelsummenabweichung des Lichstrahlendreiecks vom euklidischen
Wert proportional (und mit nahe an 1 liegendem Proportionalit\"atsfaktor) auf
die Theodolitmessungen des geod\"atischen Dreiecks $\triangle BHI$
durchschlagen würde, konnte  Gau\ss{} --- zumindest in einem  {\em
heuristischen Sinne} --- problemlos annehmen. Falls er sich die Frage genauer stellte, wa eine Absch\"atzung, wie hier im Anhang erl\"autert, ein ``neben dem Morgenkaffee'' zu erledigendes Problem.
Nimmt man die Bemerkung an Taurinus, die von Sartorius zitierte
und andere zusammen und liest sie im Kontext seiner geod\"atischen Arbeiten
der frühen 1820er Jahre, so wird man kaum 
daran zweifeln können, dass Gau\ss{} sich mit Fragen dieser Art befasste. Die Einsch\"atzung, ob Gau\ss{}  sich genau mit der {\em hier diskutierten} Frage besch\"aftigt hat, bleibt der Leserin überlassen, solange  nicht eventuell noch eigenh\"andige schriftliche Dokumente von Gau\ss{} gefunden werden, die dies eindeutig entscheiden. Ich selber  bewerte die vorliegenden Quellen samt Rekonstruktion der für ihn einsehbaren Theorieperspektive  als  einen {\em guten} historischen {\em Indizien-}Beweis.  
\vspace{0.5mm}\\
Wir können
also festhalten: \vspace{1.5mm}\\
{\em Der Inhalt von Sachverhalt  4a  bzw. 4b war für Gau\ss{} zumindest heuristisch gut
absch\"atzbar. Die zugrunde liegende Frage oder eine nahe verwandte, wird er sich gestellt haben müssen,
 wenn man seine Äu\ss{}erungen zur empirischen Bedeutung
der  Grundlagenfragen der Geometrie berücksichtigt und ernst nimmt.}\\

Wir kommen so zu der Schlussfolgerung, dass es für Gau\ss{} sehr nahe lag  zu
prüfen, ob  bei den grö\ss{}ten auf der Erde vermessbaren Dreiecken ein
systematischer Anteil des Schlie\ss{}ungsfehlers feststellbar ist, der auf eine
schon auf terrestrischem Niveau wahrnehmbare Abweichung der Raumkrümmung von
0 hinweist. {\em Die Antwort auf diese Frage war nach den Messungen von 1823 klar
negativ.} Die Schlie\ss{}ungsfehler  einiger kleiner Dreiecke des
Hannoveranischen Netzes lagen zum Teil über $1'$ und übertrafen damit den des ``gro\ss{}en'' Dreiecks  erheblich; und selbst noch der ``mittlere Fehler'' der bei den Ausgleichsrechnungen ermittelten Richtungskorrekturen  des gesamten Netzes (ohne Einbezug des ``gro\ss{}en'' Dreieicks)
lag bei $0,48''$, wie Gau\ss{}  1823 im November an Bessel
mitteilte,\footnote{\cite[IX, 366]{Gauss:Werke}, zitiert auch in \cite[282]{Breitenberger:Gauss}.}
und damit nur wenig unter dem des $\triangle BHI$.

Allerdings ist noch eine kritische Rückfrage zu bedenken: Warum gibt es
keinen Hinweis auf Versuche von  Gau\ss{}, die Winkel des Lichtstrahlendreiecks
$\triangle\hat{B} \hat{H} \hat{I} $  direkt zu messen, wenn er wirklich
Interesse daran hatte, an diesem Fall die Grö\ss{}enordnung der von ihm
vermuteten  Raumkrümmung  zu überprüfen?\footnote{Wieder so eine Frage
Jeremy Grays, die direkt den Kern unserer Geschichte erhellt.} W\"are das
nicht einfacher, als all die verwickelten Überlegungen zum Verh\"altnis von
euklidischer und nichteuklidischer Geometrie in Kauf zu nehmen, die hier zu
diskutieren waren und die die mathematikhistorische Diskussion so
irritiert haben? Und ist nicht das Fehlen jedes Nachweises einer direkten
Winkelmessung  des ebenen gro\ss{}en Dreiecks ein Hinweis darauf, dass doch
alles anders war, und Gau\ss{} nichts anderes erreichen wollte als einen besseren
Schlie\ss{}ungsfehler als alle Vorg\"anger, aber eben doch nur aus Sicht   der
traditionellen Geod\"asie?

So nahe diese Frage an den Kern der Überlegungen herankommt, so wenig
liefert sie einen weiteren (den  letzten ?)  Einwand gegen die Glaubwürdigkeit der von Sartorius berichteten Gau\ss{}schen Auswertungsüberlegung  der Messungen von 1823. Man muss sich nur vergegenw\"artigen, was  Gau\ss{} h\"atte tun müssen, um die
Winkel der Ebene $\hat{B} \hat{H} \hat{I}$ direkt  zu vermessen. Noch die
kleinste Schwierigkeit w\"are gewesen, dass er König Georg IV die
Finanzierung eines dritten Messtrupps  h\"atte schmackhaft machen müssen
(obwohl selbst das möglicherweise nicht leicht gewesen w\"are, weil ein geod\"atischer Nutzen nicht zu vesprechen war). Die
Heliotropmessungen machen dann ja nicht nur im Hauptmesspunkt (etwa auf dem
Hohehagen bei der Messung des Winkels $\angle \hat{B} \hat{H} \hat{I}$)
einen Messtrupp notwendig,  sondern je einen weiteren in den beiden anderen
Dreieckspunkten (auf dem Brocken und dem Inselsberg), um die  beiden
Heliotropen zur Markierung der Zielpunkte der Visierlinien einzurichten.
Aber das Hauptproblem w\"are ein anderes gewesen: Es w\"are n\"amlich ein völlig
neuartiges Messinstrument erforderlich geworden, mit dem simultan oder in
kurzer zeitlicher Aufeinanderfolge  zwei optische Achsen ausgerichtet und eingemessen
werden können, um dann den dazwischen liegenden  Winkel in einer geneigten
Ebene mit höchster Pr\"azision zu vermessen. Mit einem Theodoliten, selbst mit
dem von Gau\ss{} verbesserten, war und bleibt eine solche Messung
unmöglich.\footnote{An dieser Stelle wirkt sich die irreführende Darstellung
der von Gau\ss{} angeblich ausgeführten Schritte bei der Auswertung der 
Messungen  von Herrn \citeasnoun[280]{Breitenberger:Gauss} erheblich aus (siehe Anm.
(22). Sie erweckt den Eindruck, als w\"aren  die  Winkel des
Lichtstrahlendreiecks messmethodisch prim\"ar und die projizierten Winkel des
geod\"atischen Dreiecks daraus im ``zweiten Schritt'' erst nachtr\"aglich zu
berechnen. In Wirklichkeit ist jedoch  die  Projektion der Winkelmessung in die Tangentialebene durch die Horizontierung des Theodolits gewisserma\ss{}en schon  in den Messprozess eingebaut. Die (dazu schief liegenden) Winkel des Lichtstrahlendreiecks können  mit dem Theodoliten überhaupt nicht direkt gemessen werden. Ihre  Berechnung ist weder für die geod\"atische Auswertung  noch für deren  Auswertung als Genauigkeitsschranke für die Gültigkeit der euklidischen Geometrie nötig, wenn man das  Schlie\ss{}ungkriterium  wie oben erl\"autert einsetzt.}

Es bleibt also dabei: Die in Sartorius' Bericht angedeutete {\em Methode, den
Schlie\ss{}ungsfehler selbst als ein Schrankenkriterium für die Raumkrümmung zu
verwenden}, bleibt im Rahmen der vorhandenen Messinstrumente erheblich
einfacher und durch leichte
Zuspitzung der üblichen geod\"atischen Messmethoden (Auswahl und Abtrennung
eines ``gro\ss{}en'' Dreiecks vom Netz, Messung mit der höchsten erreichbaren
Pr\"azision für Einzeldreiecke)  realisierbar. Dabei sei hier noch einmal ausdrücklich darauf hingewiesen, dass dazu keine  Rückrechnung auf die einzelnen Winkel des Lichtstrahlendreiecks $\triangle \hat{B} \hat{H} \hat{I} $ erforderlich war,\footnote{Wollte man das auf der Basis von Theodolitmessungen machen, so  w\"are eine die erreichbare Pr\"azision weit überschreitende  Genauigkeit bei der  Bestimmung der sehr kleinen Neigungswinkel des Dreiecks gegen die jeweilige Tangentialebene an das Schwerefeld der Erde notwendig. Höhenwinkel sind mit guten Instrumenten ``heute'' (2003) bis auf 0,1'' messbar; sie w\"aren aber unbrauchbar  zur Bestimmung der schiefen Winkel des  $\triangle \hat{B} \hat{H} \hat{I} $, weil die atmosph\"arische Refraktion über so lange Wege sehr schwer kontrollierbar ist. Daher würde der systematische Fehler bei kleinen Höhenwinkeln das Ergebnis schlicht unbrauchbar machen (Auskunft von Herrn Heindl).} 
sondern allein die Abweichung der {\em Winkelsumme}  berücksichtigt werden musste.  Dabei war es lediglich eine Frage des Geschmackes,  der Anschaulichkeit oder der direkten Vergleichbarkeit, den Schlie\ss{}ungsfehler bezüglich der erwarteten sph\"arischen Winkelsumme ($180^{\circ} 0' 14,85\ldots ''$) oder der um den sph\"arischen Exzesss bereinigten ebenen Winkelsumme (``zwei Rechte'') anzugeben.\vspace{1mm}\\
{\em So gibt uns das Ergebnis dieser
Rücküberlegung noch einmal eine Best\"atigung dafür, dass das von Sartorius
angegebene  Kriterium auf eine höchst sinnvolle Methode der Überprüfung
hinweist. }

Dies gilt  allerdings nur, wenn man  bereit ist, die
Gau\ss{}sche Frage nach einer möglichen Grö\ss{}enordnung der Raumkrümmung in einem
``Verh\"altnissse zu solchen Grössen, welche im Bereich unserer Messungen auf
der Erde'' liegen, überhaupt als ein  sinnvolles Problem anzuerkennen.

\subsubsection*{6.  Zusammenfassung und Ausblick}
Gau\ss{} h\"atte nur dann Grund gehabt,  der Frage eines eventuellen terrestrischen
Nachweises der Krümmung des physikalischen Raumes nach 1823 weiter
nachzugehen, wenn beim gro\ss{}en Dreieck ein auff\"allig gro\ss{}er negativer
Schlie\ss{}ungsfehler aufgetreten w\"are. Tats\"achlich lag jedoch die aus  dem (negativen) Schlie\ss{}ungsfehler zu ermittelnde Korrektur der einzelnen Winkel bei $0,2''$ und damit sogar noch  unterhalb des  ``mittleren (quadratischen) Fehlers'' (der empirischen Standardabweichung $\sigma = 0, 48''$) der Ausgleichskorrekturen von Einzelrichtungen, die Gau\ss{} bei seinem Netzausgleich anzubringen hatte.\footnote{Sollte Gau\ss{} einen entsprechenden Vergleich in seinen mündlichen Mitteilungen an Sartorius von Waltershausen angestellt haben? Dies zu behaupten w\"are vielleicht zu waghalsig; es ist aber keineswegs  auszuschlie\ss{}en.  Den Wert 
$0, 48''$ als ``{\em mittlere(n) Fehler} aller Richtungen, verstanden wie in meiner {\em Theoria Combinationis}'' bezüglich der Ausgleichungen aller Hauptdreiecke der Messungen von 1821 bis 1823 gab Gau\ss{} jedenfalls schon  in seinem Briefwechsel mit Olbers am 2. 11. 1823 an \cite[II, 260]{Gauss:Olbers}. }
Aus Sicht des 20. Jahrhunderts, also hier der allgemeinen Relativit\"atstheorie,  ist dieses negative Ergebnis nicht verwunderlich. Die Absch\"atzung der Auswirkungen des  Gravitationsfeldes der Erde für die Lichstrahlengeometrie  eines  Dreiecks der Grö\ss{}enordnung des $\triangle BHI $ durch eine Schwarzschild-Metrik führt auf  eine   Winkelsumme  $\pi + 3 \cdot 10^{-13}$ \cite{Richter:Gauss}.\footnote{Den Hinweis auf diese Berechnung  verdanke ich E. Breitenberger.} 
Es ergibt sich damit eine (positive) relativistische Winkelkorrektur in der Grö\ss{}enordnung 
$10^{-8 } $$''$, also noch einmal vier Grö\ss{}enordnungen unterhalb des  von Gau\ss{} in den {\em Disquisitiones} abgesch\"atzten Nichtsph\"arizit\"atseffektes bei Anwendung des verallgemeinerten Legendre-Theorems.  

Nichts davon war zu Beginn des 19. Jahrhunderts a priori absehbar. W\"are es anders gekommen, h\"atte also Gau\ss{} eine 
 signifikante Abweichung der Winkelsumme gefunden, so w\"are natürlich weit mehr als nur  {\em eine} Auswertungsrechnung des {\em einen gro\ss{}en} Dreiecks
erforderlich geworden. Das Auftreten eines solchen von Gau\ss{}  bis  1823/24
offenbar grunds\"atzlich für möglich gehaltene Effektes w\"are wohl nur der
Auftakt für ein neues umfangreicheres Projekt geworden, die Raumkrümmung zu
vermessen. Es bleibt uns erspart zu diskutieren, welche Verbesserungen bei
der Messung der Erddaten,  der Auswahl und Koordinatenbestimmung
geeigneter Beobachtungspositionen für ``gro\ss{}e'' Dreiecke, des Ausschlusses
atmosph\"arischer Effekte, der angewendeten verbesserten Messmethoden der
Dreieckswinkel  etc. im Rahmen eines solchen fiktiven Projektes
``terrestrische Messung der Raumkrümmung'' notwendig geworden w\"aren.\footnote{Hinzu tritt natürlich, dass unter der höchst fiktiven  Annahme eines  nachweisbaren Effektes der Raumkrümmung für die Winkelsummenmessungen die Abweichung {\em positiv} gewesen w\"are.  }
 Es ist
aber ein Scheinargument, dass all dies methodologisch unmöglich gewesen w\"are; und
auch das Kopfschütteln von  besser-wisssenden Sp\"atergeborenen, ein Gau\ss{}
könne  doch nicht im irdischen Ma\ss{}stab statt im astronomischen  nach
möglichen Effekten der Raumkrümmung  gesucht haben, geht an der historischen
Situation der 1820er Jahre vorbei, in denen  glaubwürdige Parallaxdaten nicht vorlagen und nicht abzusehen war, wann sie zur Verfügung stehen würden.

Gau\ss{}' eigene (terrestrische) Messungen waren zwischen 1823 und 1838 die
verl\"asslichsten Daten zur empirischen Bestimmung einer oberen Schranke des
Betrags der Raumkrümmung (der ``Constanten'' der nichteuklidischen
Geometrie)
\[ C = |\kappa | < 10 ^{-9} \; \mbox{km}^{-2} \approx 0,04\, R^{-2} .\]
Wie oben gezeigt, ist nur eine einfache, aus Gau\ss{}' Sicht methodologisch naheliegende Dreisatzrechung erforderlich, die
von Sartorius überlieferten Gau\ss{}sche  Angabe der Schranke in diese Form der Gl. (4) zu bringen.

Auch ist gut zu verstehen, warum  sich die Gau\ss{}schen  Äu\ss{}erungen zur Grö\ss{}enordnung der von ihm erwogenen Raumkrümmung ab Mitte der 
 1820er Jahre merklich ver\"anderten. In dem nun schon h\"aufig angeführten
Brief an Taurinus im November 1824, ein Jahr nach den Messungen des Dreiecks
$BHI$, aber noch w\"ahrend der  Auswertungsarbeiten\footnote{Die Netzmessungen
liefen noch im ganzen Jahr 1824, die grundlegenden Auswertungen waren im
Januar 1825 erreicht (Gau\ss{} an Schumacher, 7.1. 1825), die
Ausgleichsrechnungen schloss Gau\ss{} erst im Mai 1826 endgültig ab (Brief an
Olbers, 14.5. 1826); siehe  \cite[100--103]{Gerardy:Gauss}.}
lie\ss{}  Gau\ss{} durchaus noch zu, wenn auch natürlich im Konjunktiv,  dass sich
``jene Constante'' gegebenenfalls durch ``Messungen auf der Erde oder am
Himmel \ldots ausmitteln'' lasse.

Am Ende seiner Auswertungen wusste er es besser. Im Sommer 1831 schrieb er
in einer --- für seine Verh\"altnisse recht ausführlichen --- Erl\"auterung über
seine Einsicht in die NEG an H.C. Schumacher von jener ``Constante'' (die
nun $k$ hie\ss{}):
\begin{quote}
$k$ [ist] eine Constante (\ldots) von der wir durch Erfahrung wissen, dass
sie gegen alles durch uns messbare ungeheuerlich gross sein muss. In Euklids
Geometrie wird sie unendlich. \cite[VIII, 215]{Gauss:Werke} zitiert nach \cite[35]{Reichardt:Gauss}
\end{quote}

Im Lichte der hier zusammengestellten Evidenz handelt es sich bei seinen unterschiedlichen
Formulierungen über die von ihm für möglich gehaltene empirische
Nachweisbarkeit der ``Constante'' nicht, wie manchmal dargestellt, um eine
Unentschiedenheit von Gau\ss{}, sondern um  Unterschiede, die in wohlformulierter
Weise eine zeitliche Entwicklung des von Gau\ss{} als gesichert angesehen Wissens
andeuten.

Dies gilt offensichtlich auch für Gau\ss{}' Hinwendung zu astronomischen
Messungen zur Bestimmung der Raumkrümmung. Ende der 1820er Jahre
war für ihn gekl\"art, dass die Konstante $| \kappa |$ 
unterhalb der   Nachweisbarkeit bei Messungen ``auf der Erde'' (lies: für
geod\"atische Messungen) liegt. Naheliegenderweise richtete sich sein
Augenmerk nun st\"arker auf astronomische Nachweismöglichkeiten.
Deren Problematik war, dass bei dem Versuch einer direkten Beobachtung die
Aberration des Lichtes (durch den Einfluss der Erdbewegung auf den
Richtungsvektor des einfallenden Strahles) eine st\"arkere scheinbare
Positionsver\"anderung hervorruft (in der Grö\ss{}enordnung $1''$) als der
eigentliche trigonometrische Parallaxeffekt (Grö\ss{}enordnung $\leq 0,1''$). Trotz
dieser schon von J. Bradley 1729 aufgewiesenen Problematik gab es zu
Beginn des 19. Jahrhunderts verschiedene Ankündigungen von
``Parallaxmessungen'' in der Grö\ss{}enordnung $1''$, die allerdings unter den Astronomen zurecht umstritten blieben. Für Gau\ss{} waren sie alle indiskutabel. Die
Lage für Parallaxmessungen \"anderte sich erst im Jahre 1838  mit F.W. Bessels
Erfolg bei der Bestimmung der parallaktischen Verschiebung  von  61 Cygni
gegenüber zwei naheliegenden anderen Sternen  ($0,314'' \pm 0,02 $). In
kurzer Zeit folgten nun weitere (stabile) Parallaxmessungen nach dieser
Methode durch W. Struve  und T. Henderson, darunter für $\alpha $ Centauri
(etwa  $0,6''$).\footnote{Siehe \cite[254, 279]{North:Astronomie}.}

Lobatschewsky hatte schon im Jahre 1830 eine Absch\"atzung der Raumkrümmung
mit einem   idealisierten, aber theoretisch grunds\"atzlich gut begründeteten
Argument vorgenommen.
Dabei berief er sich allerdings  auf  eine dieser problematischen
``Parallaxbeobachtungen'' (von $1,24''$ für den Sirius).\footnote{\cite[22f.]{Lobachevsky:Anfangsgruende}; vgl.
\cite{Daniels:Loba}.}

Wie wir wissen, begann Gau\ss{}  ab etwa 1841 mit einem Studium der  Arbeiten Lobatschewskys.
Da er ab 1838 russisch  lernte, ist es gut möglich und bei
seinem Interesse an der Problematik sogar anzunehmen, dass er in den 1840er Jahren  die entsprechende Überlegung von Lobatschewsky kennen lernte.  Unter den neuen
empirischen Bedingungen der Parallaxmessungen seit 1838 lagen nun alle
Bedingungen dafür bereit, das Problem der Raumkrümmung  durch astronomische
Messungen neu, und nun pr\"azise, zu behandeln.  Wiederum gibt es einen Zeitzeugen, diesesmal B.
Listing, der darüber berichtete, dass Gau\ss{}  in den 1840er Jahren  in seinem
Seminar über dieses Problem sprach.\footnote{Listing berichtete auch nur
mündlich in den 1870er Jahren gegenüber seinen Studenten über seine
Erinnerung daran \cite{Hoppe:Gauss}. Diese Information wurde
dankenswerterweise  von Herrn Breitenberger wieder in Erinnerung gerufen,
\cite[289]{Breitenberger:Gauss}. }
Zu diesem Zeitpunkt spielte die Erw\"ahnung eine möglichen Bestimmung der
Raumkrümmung durch terrerestrische Messungen schon  keine Rolle mehr.

Wahrscheinlich erfuhr Bernhard Riemann von dem von Gau\ss{} nun neu gestellten
Problem einer Absch\"atzung der Raumkrümmung aus Parallaxmessungen direkt oder
indirekt im Umkreis seines schwer zug\"anglichen Lehrers. Jedenfalls gab er in
seinem Habilitationsvortrag 1854 ein knappes, aber sehr pr\"agnantes Argument
für eine Absch\"atzung der Raumkrümmung aus astronomischen Messungen an. Riemann
wandte (nach einer naheliegenden Rekonstruktion)\footnote{Eine eigene Publikation zu diesem Thema ist in Vorbereitung.}
dieselbe Idee, die  Gau\ss{}  der Absch\"atzung aus terrestrischen Messungen  zugrunde gelegt hatte,  auf astronomische Beobachtungen an   --- nun allerdings theoretisch besser fundiert im  Rahmen seiner  Geometrie der (Riemannschen)  Mannigfaltigkeiten. Er erhielt so 
eine von ihm in verbaler Form, aber völlig pr\"azise angegebene Schranke für die r\"aumlich gemittelte Raumkrümmung. Mit einem etwas anderen Argument, aber mit quantitativ gleichem Ergebnis  im hyperbolischen Fall, wurde diese Schranke offenbar  unabh\"angig   von K. Schwarzschild wieder abgeleitet:  $ r > 4 \, \cdot 10^6$ Erd{\em bahn}radien (mit $r$  ``Radius''  der Raumkrümmung $r$)   \cite[345]{Schwarzschild:1900}. Die Auswertung der astronomischen Parallaxmessungen gaben damit für den Krümmungsradius Absch\"atzungen, die  die Gau\ss{}sche Schranke um 10 Grö\ss{}enordnungen überschritten  (Herabsetzung der  Gau\ss{}schen ``Constante'' $C = | \kappa |$ also um  20 Grö\ss{}enordnungen!). Aus dieser Sicht erschienen natürlich die raumgeometrischen Ergebnisse der Gau\ss{}schen Messungen der 1820er Jahre  völlig obsolet.  Doch hatte man das ja nicht im vorhinein wissen können. 

Immerhin war nun   ein neues  Kapitel der Geschichte der
empirischen Bestimmung der metrischen Raumstruktur eröffnet. Auch dieses war noch g\"anzlich 
vor-relativistisch. Allerdings waren ---  nach einer Einsch\"atzung von M. Schemmel ---  K. Schwarzschilds Betrachtungen über die NEG in der Astronomie \cite{Schwarzschild:1900} eine gute Voraussetzung für dessen rasche,  unter Astronomen exzeptionell frühe,  Annahme und produktive Bearbeitung der Relativit\"atstheorie sechzehn Jahre sp\"ater \cite{Schwarzschild:exterior,Schwarzschild:interior}. So scheint es doch  indirekte historische Beziehungen zwischen den Überlegungen  zur empirischen Bestimmung der Raumstruktur im 19. Jahrhundert und der allgemeinen Relativit\"atstheorie zu geben, wenn auch ganz anders, als dies  A. Miller  seinerzeit glaubte unterstellen zu dürfen.

\subsubsection*{Anhang: Erg\"anzung zur Begründung von  Sachverhalt 4}
Bezeichnen wir den Neigungswinkel der Ebene des Dreiecks $\triangle \hat{B} \hat{H} \hat{I} $ im Punkt  $ \hat{B}$ gegenüber der Normalen zum  Geoid mit $\nu $, so stehen der Winkel $\hat\beta $ im geneigten Dreieck und der in die Tangentialebene  projizierte Winkel $\beta $  in der Beziehung 
\[ \tan \frac{\beta }{2}   = \frac{1}{\sin \nu } \tan \frac{\hat{\beta }}{2} . \]

\begin{figure}[here]
\begin{picture}(12, 8)
 \includegraphics[height=7cm]{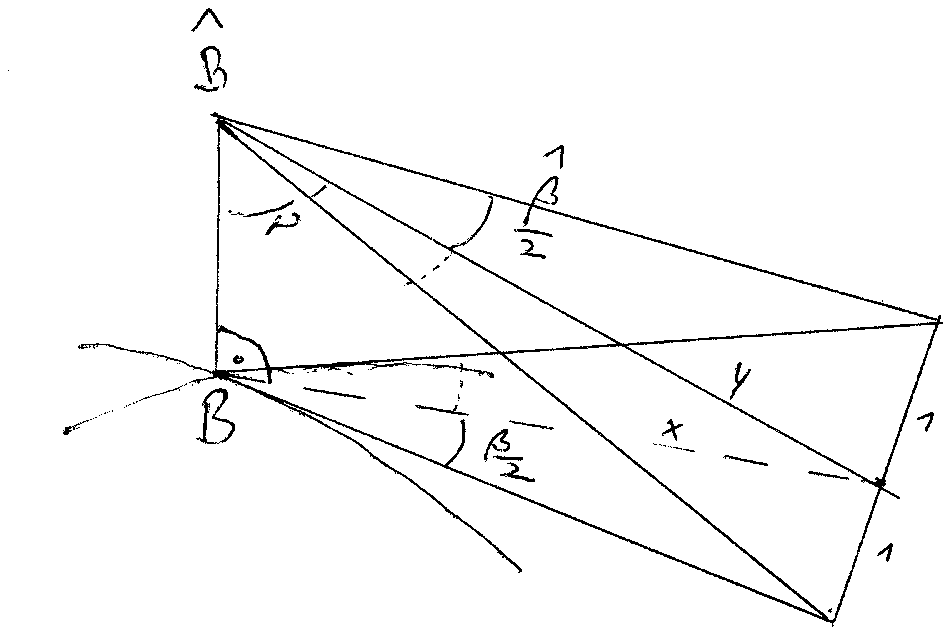} 
\end{picture} 
\caption{Winkel im geod\"atischen Dreieck und Lichtstrahlendreick}
\end{figure}

\newpage

Mit den Abkürzungen
\[  c := \frac{1}{\sin \nu } , \; \;  z_0 := \hat{\beta }_0,  \; \; h:=  \Delta \beta , \;\;z = z_0 + h = \hat{\beta }\] 
wird
\[ \beta (z) = 2 \arctan ( c \tan \frac{z}{2} ) .\]
  Reihenentwicklung  um $z_0$ liefert für $\beta _0 := 2 \arctan ( c \tan \frac{z_0}{2} ) $:
\beq
\beta (z) =  \beta _0 + \frac{c (1 + \tan ^2 \frac{z_0}{2})}{1 + c^2 tan ^2 \frac{z_0}{2}}  h + o(h)
\eeq
Beim   Gau\ss{}schen ``gro\ss{}en'' Dreieck differiert  $\nu $  um weniger als 0,5° vom rechten Winkel;  $c$ unterscheidet sich  erst in der Gr\"o\ss{}enordnung $10^{-5}$ von  1. Der Koeffizient des  linearen Term der Reihenentwicklung ist daher in sehr guter N\"aherung   $c \approx 1$.
 \bibliographystyle{apsr}
  \bibliography{a_litfile}

\subsubsection*{ }
\end{document}